\documentclass[notitlepage,leqno,10pt]{article}
\usepackage{amssymb}
\catcode`\@=11 \@addtoreset{equation}{section}

\catcode`\@=12
\usepackage[utf8]{inputenc}
\usepackage{latexsym}
\usepackage[english,francais]{babel}
\usepackage{textcomp}
\usepackage{amsmath}
\usepackage{amsfonts}
\usepackage{amssymb}
\usepackage{mathrsfs}
\renewcommand{\a }{\alpha}
\renewcommand{\b }{\beta}
\renewcommand{\d}{\delta}

\newcommand{\D }{\Delta }

\newcommand{\e }{\varepsilon }

\renewcommand{\l }{\lambda }

\newcommand{\n }{\nabla }
\newcommand{\vp }{\varphi }

\newcommand{\s }{\sigma }

\renewcommand{\O }{\Omega }

\newcommand{\be}{\begin{equation}}
\newcommand{\ee}{\end{equation}}

\newcommand{\R}{\mathbb{R}}
\newcommand{\N}{\mathbb{N}}

\newcommand{\C}{\mathbb{C}}

\newcommand{\M}{\mathcal{M}}

\newcommand{\calC }{\mathcal{C}}

\newtheorem{Theorem}{Theorem}[section]
\newtheorem{Lemma}[Theorem]{Lemma}
\newtheorem{Proposition}[Theorem]{Proposition}

\newtheorem{Remark}[Theorem]{Remark}

\def\proof{\noindent{{\bf Proof. }}}
\def\square{\vbox{
\hrule height .4pt
\hbox{\vrule width .4pt height 7pt \kern 7pt
\vrule width .4pt}
\hrule height .4pt }}
\def\QED{\hfill {$\square$}\goodbreak \medskip}
\def\R{{\mathbb R}}

\def\C{{\mathcal C}}

\font\sc=cmcsc9  \textwidth=14truecm
\author{El hadji Abdoulaye Thiam
\footnote{\footnotesize{AIMS-Senegal km2 route de Joal.
E-mail: \textit{elhadji@aims-senegal.org},
\textit{heat1719@gmail.com}. } }}
\begin{document}
\title{Hardy and Hardy-Sobolev inequalities on Riemannian manifolds}
\date{}
\maketitle
\bigskip
\noindent {\footnotesize{\bf Abstract.}
Let $ (\M,g) $ be a smooth compact Riemannian manifold of dimension $ N \geq 3 $. Given  $p_0 \in \M$, $\lambda$ $\in \R$ and $\s \in (0,2]$, we study existence and non existence of minimizers of the following quotient:
\begin{equation}\label{Paper Equation}
\mu_{\lambda,\sigma}=\inf_{u \in H^1(\M)\setminus \lbrace0\rbrace} \frac{\displaystyle\int_\M |\nabla u|^2 dv_g -\lambda \int_\M u^2 dv_g }{\biggl(\displaystyle\int_\M \rho^{-\sigma} |u|^{2^*(\sigma)} dv_g\biggl)^{2/2^*(\sigma)}},
\end{equation}
where $\rho(.):=dist(p_0,.)$ denoted the geodesic distance from $p \in \M$ to $p_0$.
In particular for $\s=2$, we provide sufficient and necessary conditions of existence of minimizers in terms of $\lambda$. For $\sigma\in (0,2)$ we prove existence of minimizers under scalar curvature pinching.\\
\textbf{Key Words}: Hardy inequality, Pure Hardy-Sobolev inequality, scalar curvature.
\section{Introduction}\label{section 1}
For $N \geq 3$, the Hardy inequality states that 
\begin{equation}\label{Hardy}
\displaystyle
\int_{\R^N} |\nabla u|^2 dx \geq \biggl(\frac{N-2}{2}\biggl)^2 \int_{\R^N} |x|^{-2} |u|^2 dx,\;\;\ \forall u \in \mathcal{D}^{1,2}({\R^N}).
\end{equation}
The constant $\biggl(\displaystyle\frac{N-2}{2}\biggl)^2 $ is sharp and never achieved in $\mathcal{D}^{1,2}({\R^N})$.  In contrast to the  Sobolev  inequality 
\begin{equation}\label{Sobolev}
\int_{\R^N} |\nabla u|^2 dx \geq S_{N,0} \biggl(\int_{\R^N} |u|^{2^*} dx\biggl)^{2/2^*}, \;\;\ \forall u \in \mathcal{D}^{1,2}({\R^N}),
\end{equation}
the best constant $$S_{N,0}= \displaystyle\frac{N(N-2)}{4}\omega_N^{2/N}$$ is achieved in $\mathcal{D}^{1,2}({\R^N})$. Here $\omega_N=|S^{N-1}|$ is the volume of the N-sphere and $ \displaystyle 2^*=\frac{2N}{N-2} $ is the critical Sobolev exponent.
For more details related to Hardy and Sobolev inequalities you can refer to the works of Brezis-Vasquez \cite{BV}, Davila-Dupaign \cite{DD1}, D'Ambrosio \cite{Amb1, Amb2}, Brezis-Marcus-Safrir \cite{BMS}, Musina \cite{Musina}, a nice exposition book in Druet-Hebey-Robert \cite{OEF} and references therein. There is also a detailed history related to Hardy inequality type in the book of Kufner-Persson \cite{KP1}.
\newline We observe that inequalities \eqref{Hardy} and \eqref{Sobolev} are scale invariant while the Sobolev inequality is additionally translation invariant. 
They have many applications in physics, spectral theory, differential geometry mathematical physics, analysis of linear and non-linear PDEs, harmonic analysis, quantum mechanic, stochastic analysis etc... For more details see the books of Lieb-Loss \cite{LL1}, Struwe \cite{Struwe} and Evans \cite{Evans}, the works of Grigor'yan \cite{Grigor}, Gkikas \cite{Gkikas}, Grigor'yan Saloff \cite{GS1} and references therein.
\newline Using H\"{o}lder's inequality, we get the interpolation between the above two inequalities called Hardy-Sobolev inequality
\begin{equation}\label{Hardy-Sobolev}
\int_{\R^N} |\nabla u|^2 dx \geq S_{N,\s} \biggl(\int_{\R^N} |x|^{-\sigma} |u|^{2^*(\sigma)} dx\biggl)^{2/2^*(\sigma)}, \;\;\ \forall u \in \mathcal{D}^{1,2}({\R^N}),
\end{equation}
where  $\sigma \in [0,2]$ and $\displaystyle 2^*(\sigma)= \frac{2(N-\sigma)}{N-2}$ is the Hardy-Sobolev critical exponent. See \cite{GK} for more details about Hardy-Sobolev inequality. We also remark that \eqref{Hardy-Sobolev} is a particular case of the Caffarelli-Kohn-Nirenberg inequality, see \cite{CKN}. The value of the best constant is
$$S_{N,\s}:= (N-2)(N-\s) \biggl[\frac{w_{N-1}}{2-\s} \frac{\Gamma^2(N-\frac{\s}{2-\s})}{\Gamma(\frac{2(N-\s)}{2-\s})}\biggl]^{\frac{2-\s}{N-\s}}$$
where $w_{N-1}$ is the volume of the $N-$sphere and $\Gamma$ is the Euler function. It was computed by Lieb \cite{LIEB} when $\s \in (0,2).$ The ground state solution is given by
\begin{equation}
\omega(x)=\biggl((N-\sigma)(N-2)\biggl)^{\frac{N-2}{2(2-\sigma)}}(1+|x|^{2-\sigma})^{\frac{2-N}{2-\sigma}}.
\end{equation}
As said previously for $\sigma=2$ the optimal constant in  \eqref{Hardy} is not attained. However there exists a "virtual ground state" $u(x)=|x|^{(2-N)/2}$ which satisfies
\begin{equation}
\Delta u+ \biggl(\frac{N-2}{2}\biggl)^2|x|^{-2} u=0 \;\;\ in \;\;\ \R^N \setminus\lbrace0\rbrace.
\end{equation}
Note also that \eqref{Hardy-Sobolev} is scale invariant. Our interest in this paper is
to study existence of minimizers of the Hardy and  Hardy-Sobolev inequalities in a Riemannian manifold. Inequality \eqref{Hardy-Sobolev} is not in general valid on a Riemannian manifold. For a compact Riemannian manifold $(\M^N,g)$, with metric $g$ and $N\geq 3$ and letting    $p_0 \in \M$, we  have
\begin{equation}
\lambda \int_\M u^2 dv_g +\int_\M |\nabla u|^2 dv_g \geq \mu \biggl(\int_\M \rho^{-\sigma} |u|^{2^*(\sigma)} dv_g\biggl)^{2/2^*(\sigma)}, \;\ \forall u\in H^1(\M),
\end{equation}
where $\rho(p)=dist_g(p,p_0)$ is the geodesic distance between $p$ and $p_0$ and $\l,\mu\in\R$ are constants depending on $\M$. The above inequality can be obtained by a simple argument of the partion of unity, see Lemma \ref{lem5} below. We then propose to study existence and non existence of minimizers of the following quotient
\begin{equation}\label{Paper inequality}
\mu_{\lambda,\sigma,p_0}=\inf_{u \in H^1(\M)\setminus\lbrace0\rbrace} \frac{\displaystyle\int_\M |\nabla u|^2 dv_g -\lambda \int_\M u^2 dv_g }{\biggl(\displaystyle\int_\M \rho^{-\sigma} |u|^{2^*(\sigma)} dv_g\biggl)^{2/2^*(\sigma)}},
\end{equation}
with $\lambda \in \R, \;\ \sigma \in (0,2]$. If there is no ambiguity, we will write $\mu_{\lambda,\sigma}$ instead of $\mu_{\lambda,\sigma,p_0}$.

In our first main result we deal with the pure Hardy problem $\s=2$. We get   the following
\begin{Theorem}\label{Theorem 1}
Let $(\M,g)$ be a smooth compact Riemannian manifold of dimension $N \geq 3$ and $p_0 \in M$. Then there exists $\l^*=\lambda^*(\M,p_0) \in \R$ such that $\mu_{\lambda,2,p_0}(\M)$ is attained if and only if $\lambda > \lambda^*(\M).$
\end{Theorem}
To explain our result and emphasize the differences between Hardy and Hardy-Sobolev inequalities in Riemannian manifolds, some definitions are in order. For an open set $\Omega \subset \M$, we put
\begin{equation}
\mu_{\lambda,\sigma}=\inf_{u \in H^1_0(\O)} \frac{\displaystyle\int_\O |\nabla u|^2 dv_g -\lambda \int_\O u^2 dv_g }{\biggl(\displaystyle\int_\O \rho^{-\sigma} |u|^{2^*(\sigma)} dv_g\biggl)^{2/2^*(\sigma)}}
\end{equation}
The existence of $\lambda^* \in \R$ is a consequence of the local Hardy:
\begin{equation}\label{LHI}
\begin{array}{cc}
\mu_{0,2}(B_g(p_0,r))=\biggl( \displaystyle\frac{N-2}{2}\biggl)^2=S_{N,2}
\end{array}
\end{equation}
which holds for  small r.
The existence and non existence of solution are based on the construction of appropriate super and sub-solution for the linear operator
$$L_\l:=-\Delta_g-\frac{(N-2)^2}{4}\rho^{-2}-\l \rho^{-2}.$$
For that we consider the geodesic normal coordinates
$$x \in B(0,r)\subset \R^N\longmapsto F(x)=\textrm{Exp}_{p_0}\bigl(\sum_{i=1}^N x^i E_i\bigl)$$
where $\textrm{Exp}_{p_0}$ is the exponential map on $\M$ and $x=(x_1,...,x_n)\in \R^N$. Using these local coordinates  we perturb the mapping $$p \longmapsto \rho^{\frac{2-N}{2}}(p)$$ to obtain $$v_a(p)=\rho^{\frac{2-N}{2}} |\textrm{log}\rho|^a$$  for $a \in \R.$ The function $|\textrm{log}\rho|^a$ allows to control the lowered order terms of the linear operator $L_\l.$ Hence a careful choice of the parameter $a$ yields super and sub-solutions to prove Theorem \ref{Theorem 1}.
However when $\s\in (0,2)$, the situation changes due to the effect of the local geometry of $\M$. Indeed we have $$\mu_{\l,\s}< S_{N,\s}$$ provided the scalar curvature is lower bounded by a function depending on the parameter $\l$.
This then allows us to prove the following
\begin{Theorem}\label{Theorem 2}
Let $(\M,g)$ be a smooth compact Riemannian manifold of dimension $N \geq 4$, $\s \in(0,2)$, $p_0 \in \M$ and $\l\in \R$ negative. We suppose that
\begin{equation}\label{MamanAissatou1}
S_g(p_0)>-6\l.
\end{equation}
Then $$\mu_{\l,\s}< S_{N,\s}$$ and it's achieved.
\end{Theorem}
We should mention that when $\sigma=0$  the above problem is related  to the well know Yamabe problem, solved by Aubin in \cite{TA}, Schoen in \cite{RS} and Trudinger in \cite{NT}. For an exposition book of such problem you can refere to the book of Druet-Hebey-Robert \cite{OEF}. \\
When $\sigma=2$, we are dealing with an eigenvalue problem for the operator $-\D_g +\mu \rho^{-2}$. A problem of   this kind was first studied    by Brezis-Marcus in \cite{BM}. See also the work of   Fall \cite{FALL}, Fall-Musina \cite{FM} and Fall-Mahmoudi \cite{FMa}. In the above mentioned paper the singularity is placed at the boundary. Hardy and Hardy-Sobolev inequalities on Riemannian manifolds have been also studied by Carron in \cite{Carron}, Adriano-Xia in \cite{AX},  Shihshu-Li in \cite{WL}, D'Ambrosio-Dipierro in \cite{DA}, E. Mitidieri  in \cite{ME} and references therein.\\
The Hardy-Soboev inequality with boundary singularities was first studied by Ghoussoub-Kang in  \cite{GK} who discovered the local  influence of the  mean curvature of the boundary in order to get a minimizer.  Further related problems, extensions and generalizations can be found in the works of Ghoussoub-Robert \cite{GR},Y. Li and Lin in \cite{LL}, Chern-Lin in \cite{CL}, Demyanov-Nazarov in \cite{DN}.\\
The paper is organized as follows: in Section  \ref{section 2} we give some preliminaries and notations, in Section \ref{section 3} we study the linear case $\sigma=2$ and in Section \ref{section 5} we study the nonlinear case $\sigma \in (0,2)$.
\section{Preliminaries  and notations}\label{section 2}
Let $p_0 \in \M$, $T_{p_0}\M$ the tangent space space of $\M$ at the point $p_0$ and by $\lbrace E_1,E_2,...,E_N\rbrace $ its standard orthonormal basis.
We consider the local parametrization of a neighborhood of $p_0$ by the usual exponential mapping
\begin{equation}\label{Exponential mapping in local coordinates}
x \ni B^N_{r_0}\mapsto F(x)=\textrm{Exp}_{p_0}\left(\sum_{i=1}^N x_iE_i\right) \in B_g(p_0,r_0)
\end{equation}
where $B^N_r \subset \R^N$ and $B_g(p_0,r) \subset \M$ denote the euclidean and geodesic balls of radii $r$. The real $r$ is suppose to be less than the injectivity radii of the manifold $\M$.
In this local parametrization, it is well known that
\begin{equation}
\rho(F(x))= |x|
\end{equation}
and 
the Laplace-Beltrami operator is given by:
\begin{equation}\label{Laplace Beltametri Operator}
\Delta_g=-g^{ij}\biggl(\frac{\partial^2}{\partial x_i \partial x_j}-\Gamma^k_{ij} \frac{\partial}{\partial x_k}\biggl)
\end{equation}
where $ \bigl \lbrace \Gamma^k_{ij} \bigl\rbrace_{1\leq i,j,k \leq N}$ are the Christoffel symbols, $ \bigl \lbrace g_{ij} \bigl\rbrace_{1\leq i,j \leq N}$ are the components of the metric g of inverse components $g^{ij}$. Moreover note that
\begin{equation}\label{l1}
\Gamma^k_{ij}(x) =O(|x|) \qquad \textrm{and} \qquad g_{ij}=\delta_{ij} -\frac{1}{3} \sum_{\a\b=1}^N R_{i\a j\b}(p_0) x_\a x_\b+O(|x|^3)
\end{equation}
where the $R_{i\a j\b}(p_0)$ denote the components of the tensor curvature at $p_0$. Then the components of the  Ricci curvature at the point $p_0$ are given by
$$
R_{ij}(p_0)=\sum_{\a=1}^N R_{i\a j\a}(p_0).
$$
The scalar curvature at $p_0$ is then
$$
S_g(p_0)=\sum_{i=1}^N R_{ii}(p_0).
$$
In the sequel, we mean by $\mathcal{D}^{1,2}(\R^N)$, the space of functions for which their gradient are square integrable in $\R^N$ and by $o(1)$ by a function depending on $n \in \N$ for which its limit at infinity is zero.
\section{Linear case: $\sigma=2$}\label{section 3}
In this section, we deal with the case $\sigma=2$. The  "virtual" ground state 
$$\omega(x)=|x|^{(2-N)/2}$$ satisfies
\begin{equation}
-\Delta \omega = \biggl(\frac{N-2}{2}\biggl)^2|x|^{-2} \omega \qquad in \qquad \R^N\setminus \lbrace0\rbrace.
\end{equation}
Using the geodesic normal coordinates, we will perturb the mapping $p\longmapsto \omega \circ F^{-1}=\rho^{\frac{2-N}{2}}(p)$ to build super-solution to get the existence of $\lambda^*.$ Moreover, with similar arguments, we will construct a subsolution which allows us to prove non existence of minimizer for $\lambda\leq \lambda^*.$
\begin{Lemma}\label{lem2}
Let
\begin{equation}\label{Perturbation}
\omega_a(x)=|x|^{\frac{2-N}{2}} |log(|x|)|^{a}
\qquad 
\textrm{ and set} \qquad v_a(F(x))=\omega_a(x).
\end{equation}
Define 
$$
L=-\Delta_g-\biggl(\frac{N-2}{2}\biggl)^2 \rho^{-2} -\lambda \rho^{-2}.
$$
Then we have 
\begin{equation}\label{SigaBeye}
Lv_a=-a(a-1)\rho^{-2} (log\rho)^{-2} v_a -\lambda \rho^{-2} v_a+O\bigl(\rho^{\frac{2-N}{2}}(-log\rho)^a\bigl) \;\ in \;\ B_g(p_0,r_0).
\end{equation}
\end{Lemma}
\proof
Let $$\varphi(t)=t^{\frac{2-N}{2}} (-log t)^a.$$
We have
\begin{equation*}\label{ppp}
\Delta \omega_a(x)=\varphi"(|x|)+\frac{N-1}{|x|} \varphi'(|x|)= -\biggl(\frac{N-2}{2}\biggl)^2\omega_a(x)|x|^{-2}+a(a-1)\omega_a(x)|x|^{-2}(log|x|)^{-2}.
\end{equation*}
Then by \eqref{Laplace Beltametri Operator} and \eqref{l1}, we obtain
\begin{equation*}\label{Ref}
\Delta_g v_a=-\Delta_{\R^N} \omega_a+O_{ij}(|x|^2) \partial^2_{ij} \omega_a+ O_k(|x|) \frac{\partial \omega_a}{\partial x_k}.
\end{equation*}
Therefore
\begin{equation*}
Lv_a=-a(a-1)\rho^{-2} (log\rho)^{-2} v_a -\lambda \rho^{-2} v_a+O\bigl(\rho^{\frac{2-N}{2}}(-log\rho)^a\bigl).
\end{equation*}
This ends the proof of the lemma. \QED
\begin{Lemma}\label{lem3}
Let $(\M,g)$ be a smooth compact Riemannian manifold of dimension $ N\geq3 $ and $p_0 \in \M$. Then there exists $r_0>0$ such that:
\begin{equation}\label{Reff1}
\int_{B_g(p_0,r_0)} |\nabla u|^2 dv_g \geq \biggl(\frac{N-2}{2}\biggl)^2 \int_{B_g(p_0,r_0)} \rho^{-2} |u|^2 dv_g+\int_{B_g(p_0,r_0)} \rho^{-2}(log \rho)^{-2} u^2 dv_g,\quad \textrm{ for all } u \in H^1_0(B_g(p_0,r_0)).
\end{equation}
\end{Lemma}
\proof
We choose $a=-1$ in \eqref{Perturbation} and let $v_{-1}=V$. We have for $r_0$ small enough, that
\begin{equation}\label{Maman}
\frac{-\Delta V}{V}\geq \biggl(\frac{N-2}{2}\biggl)^2 \rho^{-2}+\rho^{-2} (log \rho)^{-2} \quad in \quad B_g(p_0,r_0).
\end{equation}
Let $u \in \mathcal{C}^{\infty}_c (B_g(p_0,r_0))$ and consider $ \psi=\displaystyle\frac{u}{V}. $
We have
\begin{equation*}
|\nabla u|^2=|V \nabla \psi|^2+\nabla V. \nabla (V \psi^2).
\end{equation*}
By integration by parts, we have
\begin{equation}\label{Maman1}
\int_{B_g(p_0,r_0)} |\nabla u|^2 dv_g=\int_{B_g(p_0,r_0)}|V \nabla \psi|^2dv_g-\int_{B_g(p_0,r_0)} \frac{\Delta V}{V} u^2 dv_g.
\end{equation}
Therefore by \eqref{Maman} and \eqref{Maman1} we obtain
$$
\int_{B_g(p_0,r_0)} |\nabla u|^2 dv_g \geq \biggl(\frac{N-2}{2}\biggl)^2 \int_{B_g(p_0,r_0)} \rho^{-2} u^2 dv_g +\int_{B_g(p_0,r_0)} \rho^{-2} (log \rho)^{-2} u^2 dv_g \qquad \forall u \in \mathcal{C}^{\infty}_c (B_g(p_0,r_0)).
$$
Note that $\mathcal{C}^{\infty}_c(B_g(p_0,r_0))$ is dense in $H^1_0(B_g(p_0,r_0))$. This ends the proof of the lemma. \QED
\subsection{Existence of $\lambda^*$}
\begin{Proposition}\label{Proposition 1}
Let $ (\M,g) $ be a smooth compact Riemannian manifold of dimension $ N \geq 3$. Then there exists $ \lambda^* \in \R $ such that 
$$\mu_{\lambda^*,2}(\M)=\biggl(\displaystyle\frac{N-2}{2}\biggl)^2 \qquad \textrm{and} \qquad
\mu_{\lambda,2}(\M)< \biggl(\displaystyle\frac{N-2}{2}\biggl)^2 \qquad \textrm{ for all } \lambda> \lambda^*.$$
\end{Proposition}
\proof
We Claim that for all $\l \in \R$
\begin{equation}\label{Islam} 
\mu_{\lambda,2}(\M)\leq \biggl(\frac{N-2}{2}\biggl)^2
\end{equation}
Indeed, recall that, the best constant of Hardy is given by
$$\displaystyle \mu_{0,2}(\mathbb{R}^N)=\biggl(\frac{N-2}{2}\biggl)^2.$$
Then for any $ \delta >0  $ , we can find $ u_{\delta} \in \mathcal{C}_c^{\infty}(\mathbb{R}^N) $ such that
\begin{equation}\label{m}
\int_{\mathbb{R}^N} |\nabla u_{\delta}|^2 \mathrm{d} y \leq \biggl(\biggl(\frac{N-2}{2}\biggl)^2+\delta\biggl) \int_{\mathbb{R}^N} |y|^{-2} u_{\delta}^2 \mathrm{d}y.$$
Let $$ p=F\left(\e y\right) \qquad \textrm{and} \qquad v_\delta(p)=\epsilon^{\frac{2-N}{2}} u_{\delta}(\epsilon^{-1} F^{-1}(p)).
\end{equation}
For $ \epsilon $ small enough, $ v_\delta \in \mathcal{C}_c^{\infty}(\M). $
By applying the change of variable formula and using \eqref{l1} we obtain
\begin{equation}\label{8}
\mu_{\lambda,2}(\M) \leq  \frac{\displaystyle \int_{\M} {|\nabla v_\delta|}_{g}^{2} dv_g -\lambda \int_{\M}|v_\delta|_{g}^{2} dv_g}{\displaystyle\int_{\M} \rho^{-2}|v_\delta|^{2} dv_g} \leq (1+c \epsilon) \frac{\displaystyle \int_{\mathbb{R}^N} |\nabla u_{\delta}|^2 \mathrm{d} y}{\displaystyle \int_{\mathbb{R}^N} |y|^{-2} u_{\delta}^2 \mathrm{d}y}+c\epsilon^2 |\lambda|.
\end{equation}
Hence
$$\mu_{\lambda,2}(\M) \leq \bigl(1+c\epsilon\bigl)\biggl(\bigl(\frac{N-2}{2}\bigl)^2+\delta\biggl)+c\epsilon^2|\lambda|.$$
As $ \epsilon, \delta \longrightarrow 0 $ respectively, the claim \eqref{Islam} follows.\\
We claim that there exists $ \overline{\lambda} $ such that
\begin{equation}\label{10}
\begin{array}{cc}
\displaystyle \biggl( \frac{N-2}{2}\biggl)^2 \int_{\M} \rho^2 |u|^2dv_g \leq \int_{\M} |\nabla u|^2dv_g-\overline{\lambda} \int_{\M} |u|^2 dv_g, \quad\forall u \in H^1(\M).
\end{array}
\end{equation}
Indeed: for $r_0>0$ small enough, we let $\vp \in \calC^\infty_c\left(B_g(p_0, 2r_0)\right)$ such that
$$
0\leq \vp \leq 1 \qquad \textrm{ and } \qquad \vp\equiv 1 \qquad\textrm{in } \quad B_g(p_0,r_0).
 $$
Then, we have $u \varphi \in H^1_0(B_g(p_0,2r_0))$ and 
\begin{eqnarray*}
\displaystyle \int_{\M} |u|^2 \rho^{-2} dv_g
&=&\displaystyle \int_\M |u \varphi+(1-\varphi)u|^2 \rho^{-2} dv_g \\\
&=&\displaystyle \int_{\M} |u \varphi|^2 \rho ^{-2} dv_g+ \int_{\M} |\bigl(1-\varphi \bigl)u|^2 \rho ^{-2} dv_g+2\int_{\M} |u \varphi(1-\varphi)|^2 \rho ^{-2} dv_g\\\
&\leq &\displaystyle \int_{B_g(p_0,2 ro)} |u \varphi|^2 \rho ^{-2} dv_g+3\int_{B_g^C(p_0,ro)} |(1-\varphi)u|^2 \rho ^{-2} dv_g \\\
&\leq & \biggl(\frac{N-2}{2}\biggl)^{-2}\int_{B_g(p_0,2 r)} |\nabla(u \varphi)|^2 dv_g+3 \int_{B_g^C(p_0,ro)} |u|^2 \rho ^{-2} dv_g.
\end{eqnarray*}
Futhermore we have
\begin{eqnarray*}
\displaystyle \int_{B_g(p_0,2r_0)} |\nabla (u \varphi)|^2 dv_g &\leq &\int_{\M} |\varphi \nabla u+u \nabla\varphi |^2 dv_g\\
&\leqslant & \int_{\M} |\nabla u |^2 dv_g + C(p_0,r_0,N)\int_{\M} | u |^2 dv_g+\frac{1}{2}\int_{\M} |\nabla u |^2|\nabla \varphi |^2 dv_g\\
&\leqslant &\int_{\M} |\nabla u |^2 dv_g + C(p_0,r_0,N)\int_{\M} | u |^2 dv_g -\frac{1}{2}\int_{\M} |\Delta u |^2 dv_g\\
&\leqslant & \int_{\M} |\nabla u |^2 dv_g+ C(p_0,r_0,N)\int_{\M} | u |^2 dv_g.
\end{eqnarray*}
Therefore
\begin{equation}
\displaystyle \biggl(\frac{N-2}{2}\biggl)^2 \int_{\M} |u|^2 \rho^{-2} dv_g \leqslant \int_{\M} |\nabla u |^2 dv_g+C(p_0,r_0,N)\int_{\M} | u |^2 dv_g .
\end{equation}
Hence there exist $ \overline{\lambda} $ such that $$ \mu_{\bar{\lambda},2}(\M)\geqslant \biggl(\frac{N-2}{2}\biggl)^2 .$$
Since the function  $ \lambda \longmapsto \mu_\lambda $ is decreasing, we can define $ \lambda^*$ as $$\lambda^*= \sup\biggl\lbrace  \lambda \in \mathbb{R} : \mu_{\lambda,2}(\M)=\biggl(\displaystyle \frac{N-2}{2}\biggl)^2 \biggl\rbrace.$$  Then the claim follows. This ends the proof.\QED
\begin{Proposition}\label{Proposition 2}
Let $ \mathcal{M} $ be a smooth compact manifold of dimension $N\geq3$. Then 
\begin{equation}\label{Linear Case}
\mu_{\lambda,2}(\M)=\inf_{u \in H^1(\M)\setminus\lbrace0\rbrace} \frac{\displaystyle\int_\M |\nabla u|^2 dv_g -\lambda \displaystyle\int_\M u^2 dv_g }{\displaystyle\int_\M \rho^{-2} u^2 dv_g}
\end{equation}
is achieved for every $\lambda > \lambda^*$.
\end{Proposition}
\proof
Let $\lbrace u_n \rbrace_{n\geq0}$ be a minimizing sequence of \eqref{Linear Case} normalized so that:
\begin{equation}\label{p}
\int_{\M} \frac{u^2_n}{\rho^2} dv_g=1 \qquad \textrm{and } \qquad
\mu_{\lambda,2}(\M)=\displaystyle \int_{\M} |\nabla u_n|^2 dv_g-\lambda \displaystyle\int_{\M} u_n^2 dv_g+o(1).
\end{equation}
Thus $\displaystyle \lbrace u_n \rbrace_{n\geq0}$ is bounded in $ H^1(\M) $. After passing to a subsequence, we assume that there exists $u \in H^1(\M)$ such that
\begin{equation}\label{12}
\begin{array}{cc}
u_n \rightharpoonup u \;\;\ weakly \;\ in \;\  H^1(\M).
\end{array}
\end{equation}
Let  
\begin{equation}\label{q} v_n=u_n-u .
\end{equation}
Then we have:
\begin{equation}\label{13}
v_n\longrightarrow 0 \;\;\ in \;\;\ L^2(\M)\;\;\ , \;\;\ v_n\rightharpoonup 0\;\;\ in \;\;\ H ^1(\M)\;\;\ and\;\;\   \frac{v_n}{\rho} \rightharpoonup 0\;\;\ in\;\;\  L^2(\M).
\end{equation}
Using \eqref{p} and the weak convergence we have
\begin{equation}\label{14}
 \mu_{\lambda,2}+o(1)= \displaystyle \int_{\M} |\nabla u_n|^2 dv_g-\lambda \int_{\M} u_n^2 dv_g=\int_{\M} |\nabla u|^2 dv_g+\int_{\M} |\nabla v_n|^2 dv_g-\lambda \int_{\M} u^2 dv_g+o(1).
 \end{equation}
and  that 
\begin{equation}\label{15}
1=\displaystyle \int_{\M}\frac{ u_n^2}{\rho^2} dv_g=\int_{\M} \frac{u^2}{\rho^2} dv_g+\int_{\M} \frac{v_n^2}{\rho^2} dv_g+o(1).
\end{equation}
By Proposition \ref{Proposition 1} and  \eqref{p}  we obtain
\begin{equation}\label{16}
\displaystyle \int_{\M} |\nabla v_n|^2 dv_g +o(1)\geq \left(\frac{N-2}{2}\right)^2\biggl(\int_{\M} \frac{v_n^2}{\rho^2} dv_g\biggl)=\left(\frac{N-2}{2}\right)^2\biggl(1-\int_{\M} \frac{u^2}{\rho^2} dv_g\biggl)+o(1).
\end{equation}
By \eqref{16} and  \eqref{14} we obtain:
\begin{equation}\label{17}
 \displaystyle \int_{\M} |\nabla u|^2 dv_g+ \biggl(\frac{N-2}{2}\biggl)^2\biggl( 1-\int_{\M} \frac{u^2}{\rho^2} dv_g\biggl)-\lambda \int_{\M} u^2 dv_g \leqslant \mu_{\lambda,2}.
 \end{equation}
Note that $$\displaystyle \int_{\M} |\nabla u|^2 dv_g-\lambda \int_{\M} u^2 dv_g \geq \mu_{\lambda,2} \int_{\M} \frac{u^2}{\rho^2} dv_g.$$
Then
$$\biggl(\mu_{\lambda,2}-(\frac{N-2}{2})^2\biggl)\biggl(\int_{\M} \frac{u^2}{\rho^2} dv_g-1\biggl) \leq0 .$$
Since $\displaystyle  \mu_{\lambda,2} <\biggl(\frac{N-2}{2}\biggl)^2 $ (See Proposition \ref{Proposition 1} above) we have $$ 1 \leq\displaystyle \int_{\M} \frac{u^2}{\rho^2} dv_g.$$ 
Therefore  $$\displaystyle  \int_{\M} \frac{u^2}{\rho^2}dv_g= 1.$$
So $u$ is a minimizer for $\displaystyle\mu_{\lambda,2}(\M)$ and $$\displaystyle\int_{\M}|\nabla v_n|^2 dv_g \rightarrow 0.$$  Thus $u_n \rightarrow u $ in $ H^1(\M)$ and the proof is complete. \QED
\begin{Proposition}\label{Proposition 3}
Let $ \mathcal{M} $ be a smooth compact manifold of dimension $N\geq 3$. Then $ \mu_{\lambda,2}(\M) $ is not achieved for every $\lambda \leq \lambda^*.$
\end{Proposition}
\proof
We study separately the case $ \lambda=\lambda^* $ and the case $ \lambda< \lambda^* .$
For every $ \lambda < \lambda^* $ the statement is verified.
Indeed suppose that for some $ \bar{\lambda} < \lambda^* $ the infimum is attained by $ \bar{u} \in H^1(\M) $.
We suppose that $ \bar{u} $ is normalized so that:
$$
\displaystyle \int_{\M} \frac{\bar{u}^2}{\rho^2} dv_g=1
\qquad
\textrm{ and } \qquad
  \int_{\M} |\nabla \bar{u}|^2 dv_g-\bar{\lambda} \displaystyle\int_{\M}\bar{u}^2 dv_g=\biggl(\frac{N-2}{2}\biggl)^2.
$$
Then, for $ \bar{\lambda} < \lambda < \lambda^*$ we have,
$$ 
\begin{array}{ll}
\displaystyle\biggl(\frac{N-2}{2}\biggl)^2=\mu_\lambda \leq  \displaystyle\int_{\M} |\nabla \bar{u}|^2 dv_g-{\lambda} \int_{\M}\bar{u}^2  dv_g< \biggl(\frac{N-2}{2}\biggl)^2.
\end{array}
$$
So for $ \lambda < \lambda^* $ we have $ \mu_{\lambda,2} $ is not achieved.\\
We suppose by contradiction that for $ \lambda=\lambda^* $, there exists $ u \in H^1(\M) $ such that $ \mu_{\lambda^*,2} $ is achieved.
Recall that for $ u \in H^1(\M) $ , $ |u| \in H^1(\M) $ and  $ |\nabla u|=|\nabla |u|| $ almost everywhere, see \cite{OEF}.
So we may assume that  $ u>0 $.
Let
\begin{equation}
L:=-\D_g-\biggl( \frac{N-2}{2}\biggl)^2\rho^{-2} -\lambda \rho^{-2}.
\end{equation}
By standard regularity theory, see \cite{GT} and thanks to the maximun principle
u is smooth and positive in $ \M \setminus\lbrace p_0 \rbrace$.
Recall from Lemma \ref{lem2} that
\begin{equation}\label{ll}
Lv_a=\displaystyle -a(a-1) \rho^{-2}(-log \rho)^{-2} v_a-\lambda \rho^{-2} v_a+ O\biggl(\rho^{\frac{2-N}{2}}(log \rho)^{a}\biggl).
\end{equation}
The dominant term in the right hand side of equation \eqref{ll} is  $-a(a-1) \rho^{-2}(-log \rho)^{-2} v_a$.
So for r small enough, we have for $a<-\frac{1}{2}$
$$ Lv_a\leq 0 \;\;\;\;\ in \;\;\;\;\;\ B(p_0,r)$$
and also $v_a \in H^1(B(p_0,r)).$
Now let $\epsilon >0$ such that :
$$ \epsilon v_a \leq u \;\;\;\;\;\ on \;\;\;\;\;\ \displaystyle \Sigma_r=\bigl \lbrace p \in \M : \rho(p)=r \bigl\rbrace$$
and let $$ W_a=\epsilon v_a -u.$$
Then $\displaystyle W_a^+\in H^1_0(B(p_0,r)) \;\;\;\;\ for \;\;\; all \;\;\;\ a \in \bigl(-1,-\frac{1}{2} \bigl).$
Furthermore  $$Lu\geq 0.$$
 Therefore $$ LW_a\leq 0 \;\;\;\;\ in \;\;\;\;\;\ B(p_0,r) \;\;\;\;\ \forall \;\;\;\;\ a \in \biggl(-1,- \frac{1}{2}\biggl).$$
Using \eqref{Reff1}, we deduce that
\begin{equation}\label{ii}
\int_{B(p_0,r)} \biggl( |\nabla W_a^+|^2-\bigl( \frac{N-2}{2}\bigl)^2 \rho^{-2}(W_a^+)^2    \biggl) \geq 0.
\end{equation} 
Therefore, the fact that 
$$
\begin{array}{ll}
\displaystyle \int_{B_g(p_0,r)} \biggl( |\nabla W_a^+|^2-\bigl( \frac{N-2}{2}\bigl)^2 \rho^{-2}(W_a^+)^2-\lambda \rho^{-2} (W_a^+)    \biggl)\leq 0
\end{array}
$$
implies $$\epsilon v_a \leq u  \;\;\;\;\;\;\ in \;\;\;\;\;\ B_g(p_0,r).$$
Hence $$ \epsilon \bigl( \rho^{\frac{2-N}{2}} (log \rho)^{-1} \bigl)^{\frac{1}{2}}\leq u \;\;\;\;\;\ in \;\;\;\;\;\ B_g(p_0,r) $$ and consequently $$ \displaystyle \frac{u}{\rho} \notin L^2(B(p_0,r)) .$$
This contradicts the assumption that $ u \in H^1(\M) $. \QED
\subsection{Proof of Theorem 1.1}\label{section 4}
The existence of $\lambda^*$ is given by the Proposition \ref{Proposition 1}. The proof of the "if" part is done in Proposition \ref{Proposition 2} and the "only if" part is done in Proposition \ref{Proposition 3}.
\section{Nonlinear case: $\sigma \in (0,2)$}\label{section 5}
We recall the Hardy-Sobolev best constant on the euclidean space 
\begin{equation}\label{Sensigma}
S_{N,\sigma}=\inf_{u \in \mathcal{D}^{1,2}(\R^N)} \frac{\displaystyle\int_{\R^N} |\nabla u|^2 dx}{\displaystyle\biggl(\int_{\R^N} |x|^{-\sigma} |u|^{2^*(\sigma)} dx\biggl)^{2/2^*(\sigma)} }.
\end{equation}
We will need the following
\begin{Lemma}\label{lem5}
Let $ (\M,g) $ be a smooth compact Riemannian manifold of dimension $N\geq 3$. 
For all $ \epsilon >0 $ small, there exist $ K(\epsilon, \M) $ positive constants such that for all $ u \in H^1(\M)$,
\begin{equation}\label{a}
S_{N,\sigma}\biggl( \int_{\M} \rho^{-\sigma}|u|^{2^*(\sigma) }dv_g\biggl)^{2/2^*(\sigma)} \leq (1+\epsilon) \int_\M |\nabla u|^2 dv_g+K(\epsilon,\M)\biggl[ \int_\M |u|^2 dv_g+
\biggl( \int_{\M} |u|^{2^*(\sigma) }dv_g\biggl)^{2/2^*(\sigma)}\biggl].
\end{equation}
\end{Lemma}
\proof
Let $\epsilon >0$ and $ \varphi \in \calC^\infty_c(B_g(p_0, 2\e)) $ such that $$ 0\leq \varphi \leq 1  \qquad
 \textrm{ and } \qquad \vp\equiv 1 \quad\textrm{ in } B_g(p_0,\e).
$$
We have for 
$ 2^*(\sigma) > 1,$  there exists $ C(\epsilon)>0 $ such that
\begin{equation}\label{CCC1}
|u|^{2^*(\s)}=|u\varphi+(1-\varphi)u|^{2^*(\sigma)} \leq (1+\epsilon) |u\varphi|^{2^*(\sigma)} +C(\epsilon) |(1-\varphi)u|^{2^*(\sigma)}.
\end{equation}
Then
\begin{equation}\label{DDD1}
\biggl(\displaystyle \int_\M |u|^{2^*(\sigma)} \rho^{-\sigma} dv_g\biggl)^{2/{2^*(\sigma)}} \leq (1+\epsilon) \biggl(\int_{B_g(p_0,2\epsilon)} |u \varphi|^{2^*(\sigma)} \rho^{-\sigma} dv_g\biggl)^{2/{2^*(\sigma)}} + C(\epsilon) \biggl( \int_\M |(1-\varphi)u|^{2^*(\sigma)} dv_g\biggl)^{2/{2^*(\sigma)}}.
\end{equation}
By change of variable formula, there exists a constant $C>0$ such that 
\begin{equation}\label{EEE1}
\biggl(\int_{B_g(p_0,2\epsilon)} |(u \varphi)(p)|^{2^*(\sigma)} \rho^{-\sigma}(p) dv_g \biggl)^{2/2^*(\sigma)}\leq (1+C\epsilon)\biggl(\int_{B^N_{2\epsilon}} |(u \varphi)(F(x))|^{2^*(\sigma)} |x|^{-\sigma} dx \biggl)^{2/2^*(\sigma)}
\end{equation}
and by \eqref{Sensigma}, we have that
\begin{equation}\label{FFF1}
S_{N,\s}\biggl(\int_{B(p_0,2\epsilon)} |(u \varphi)(p)|^{2^*(\sigma)} \rho^{-\sigma}(p) dv_g \biggl)^{2/2^*(\sigma)}\leq(1+C\epsilon) \int_{\R^N} |\nabla (u \varphi)(F(x))|^2 dx.
\end{equation}
Since
\begin{equation}\label{GGG1}
|\nabla (u \varphi)|^2 =|\varphi \nabla u|^2+|u \nabla \varphi|^2+ 2 u \varphi \nabla u \nabla \varphi
\end{equation}
we have
\begin{equation}\label{HHH1}
\int_{\R^N} |\nabla(u \varphi)(F(x))|^2 dx \leq \int_{B_g(p_0,2\epsilon)} |\nabla u|^2 dv_g+ C'(\epsilon,\M) \int_{B_g(p_0,2\epsilon)} |u|^2 dv_g,
\end{equation}
Hence using \eqref{DDD1}, \eqref{EEE1}, \eqref{FFF1} and \eqref{HHH1}, we get the result
\begin{equation}
S_{N,\sigma}\biggl( \int_{\M} \rho^{-\sigma}|u|^{2^*(\sigma) }dv_g\biggl)^{2/2^*(\sigma)} \leq (1+\epsilon) \int_\M |\nabla u|^2 dv_g+K(\epsilon)\biggl[ \int_\M |u|^2 dv_g+
\biggl( \int_{\M} |u|^{2^*(\sigma) }dv_g\biggl)^{2/2^*(\sigma)}\biggl],
\end{equation}
where $K(\epsilon)= Max(C'(\epsilon), C(\epsilon))$. This ends the proof. \QED
\begin{Remark}
For all $u \in \C^1(\M)$,  there exists a constant $C(\M,N)$ such that
\begin{equation}
C(\M,N) \biggl(\int_\M \rho^{-\s} |u|^{2^*(\s)} dv_g\biggl)^{2/2^*} \leq \int_\M |\nabla u|^2 dv_g+\int_\M u^2 dv_g.
\end{equation}
Indeed, using the fact that there exists a constant $K(\M,N)$ such that
\begin{equation}
K(\M,N) \biggl( \int_\M |u|^{2^*(\s)} dv_g\biggl)^{2/2^*(\s)} \leq \int_\M |\nabla u|^2 dv_g+\int_\M u^2 dv_g 
\end{equation}
and inequality \eqref{a}, the remark follows.
\end{Remark}
In particular $\mu_{\l,\s}$ is well defined for all $\lambda<0$.
\subsection{Existence Result}}
Recall that 
\begin{equation*}
\mu_{\lambda,\sigma}=\inf_{u \in H^1(\M)\setminus \lbrace 0\rbrace} \frac{\displaystyle\int_M |\nabla u|^2 dv_g -\lambda \int_M u^2 dv_g }{\biggl(\displaystyle\int_\M \rho^{-\sigma} |u|^{2^*(\sigma)} dv_g\biggl)^{2/2^*(\sigma)}}.
\end{equation*}
Note that the proof of Theorem \ref{Theorem 2} is a direct consequence of Proposition \ref{Proposition 4} and Proposition \ref{Proposition 5} below. Then we have
\begin{Proposition}\label{Proposition 4}
Let $(\M,g)$ be a smooth compact Riemannian manifold of dimension $ N \geq3$.
If $ \mu_{\lambda,\sigma} < S_{N,\sigma}$ then $ \mu_{\lambda,\sigma}$ is attained.
\end{Proposition}
\proof
Let $\lbrace u_n\rbrace_{n\geq0} $ be a minimizing sequence normalized so that
\begin{equation}\label{Normalization}
\int_\M \rho^{-\sigma} u_n^{2^*(\sigma)} dv_g=1 \qquad \textrm{ and } \qquad
 \mu_{\l,\s}=\int_\M |\n u_n|^2 dv_g-\l \int_\M u_n^2 dv_g+o(1).
\end{equation}
Then $ \lbrace u_n \rbrace_{n \geq 0}$ is bounded in $ H^1(\M)$ and we assume that up to a subsequence
\begin{equation}\label{Convergence}
u_n \rightharpoonup u \quad in\quad H^1(\M) \quad\textrm{ and }\quad u_n \longrightarrow u \quad in \quad L^{2^*(\sigma)}(\M) \quad \textrm{ for } \quad 0 < \sigma \leq 2.
\end{equation}
By the convergence in \eqref{Convergence} and the normalization \eqref{Normalization} we have
\begin{equation}\label{Estimation 1}
\mu_{\lambda, \sigma} +o(1)=\int_\M |\nabla u_n|^2 dv_g-\lambda \int_\M u^2 dv_g=\int_\M |\nabla u|^2_g dv_g+ \int_\M |\nabla (u_n-u)|^2_g dv_g -\lambda \int_\M u^2 dv_g +o(1).
\end{equation}
By Brezis-Lieb Lemma
\begin{equation}\label{Estimation 2}
1=\int_\M \rho^{-\sigma} |u_n|^{2^*(\sigma)} dv_g= \int_\M \rho^{-\sigma} |u|^{2^*(\sigma)} dv_g+\int_\M \rho^{-\sigma} |u_n-u|^{2^*(\sigma)}  dv_g+o(1).
\end{equation}
From lemma \ref{lem5} and \eqref{Convergence}, we obtain
\begin{equation}\label{Estimation 3}
S_{N,\sigma}\biggl(\int_\M \rho^{-\sigma} |u_n-u|^{2^*(\sigma)} dv_g\biggl)^{2/2^*(\sigma)} \leq (1+\epsilon) ||\nabla (u_n-u)||^2_2+o(1).
\end{equation}
Therefore
\begin{equation}
S_{N,\sigma} \biggl(1- \int_\M \rho^{-\sigma} |u|^{2^*(\sigma)} dv_g\biggl)^{2/2^*(\sigma)}\leq (1+\epsilon) ||\nabla (u_n-u)||^2_2+o(1).
\end{equation} 
From \eqref{Estimation 1} and \eqref{Estimation 3}, we get
\begin{equation}
\int_\M |\nabla u|^2 dv_g+ \frac{S_{N,\sigma}}{1+\epsilon}\biggl(1-\int_\M \rho^{-\sigma} |u|^{2^*(\sigma)} dv_g\biggl)^{2/2^*(\s)}-\lambda \int_\M u^2 dv_g \leq \mu_{\lambda,\sigma}.
\end{equation}
Since
\begin{equation}
\mu_{\lambda,\sigma} \biggl( \int_\M\rho^{-\sigma} |u|^{2^*(\sigma)} dv_g\biggl)^{2/2^*(\sigma)} \leq \int_\M |\nabla u|^2 dv_g-\lambda \int_\M u^2 dv_g ,
\end{equation}
we get
\begin{equation}
\frac{S_{N,\sigma}}{1+\epsilon}\biggl(1-\int_\M \rho^{-\sigma} |u|^{2^*(\sigma)} dv_g\biggl)^{2/2^*(\s)} \leq \mu_{\l,\s} \biggl(1-\biggl( \int_\M\rho^{-\sigma} |u|^{2^*(\sigma)} dv_g\biggl)^{2/2^*(\sigma)}\biggl).
\end{equation}
Moreover $$\displaystyle 1-\biggl( \int_\M\rho^{-\sigma} |u|^{2^*(\sigma)} dv_g\biggl)^{2/2^*(\sigma)} \leq \biggl(1-\int_\M \rho^{-\sigma} |u|^{2^*(\sigma)} dv_g\biggl)^{2/2^*(\s)}.$$
Taking the limit as $\epsilon \longrightarrow 0$ we obtain
\begin{equation}
\bigl( S_{N,\s}- \mu_{\l,\s}\bigl) \biggl(1-\biggl( \int_\M\rho^{-\sigma} |u|^{2^*(\sigma)} dv_g\biggl)^{2/2^*(\sigma)}\biggl) \leq 0.
\end{equation}
Since $$S_{N,\s}< \mu_{\l,\s}\qquad \textrm{ and } \qquad \displaystyle \int_\M\rho^{-\sigma} |u|^{2^*(\sigma)} dv_g \leq 1$$
it follows $$ \displaystyle \int_\M \rho^{-\sigma} |u|^{2^*(\sigma)}= 1.$$
Therefore $u_n \longrightarrow u$ in $H^1(\M)$. In particular $u$ is a minimizer for $\mu_{\lambda,\sigma}$. \QED
In the following we give necessary condition to get strict inequality between Hardy-Sobolev best constant $\mu_{\lambda,\sigma}$ and $ S_{N,\s}$ in order to get sufficient condition for existence of minimizer. Then we have
\begin{Proposition}\label{Proposition 5}
Let $(\M,g)$ be a smooth compact Riemannian manifold of dimension $N \geq 4$, $\s \in(0,2)$, $p_0 \in \M$ and $\l\in \R$ negative. We assume that
\begin{equation}\label{MamanAissatou}
S_g(p_0)>-6\l.
\end{equation}
Then
\begin{equation}\label{aa}
\mu_{\lambda,\sigma,p_0}=\mu_{\lambda,\sigma} < S_{N,\sigma}.
\end{equation}
\end{Proposition}
\proof
Let
\begin{equation}\label{GroundState}
w(x)=\bigl( 1+|x|^{2-\sigma}\bigl)^{\frac{2-N}{2-\sigma}}
\end{equation}
the  ground-state solution of the best Hardy-Sobolev constant
\begin{equation}\label{Hardy-SobolevBestConstant}
S_{N,\sigma}=\frac{\displaystyle \int_{\R^N} |\nabla w|^2 dx}{\displaystyle \biggl(\int_{\R^N} |x|^{-\sigma} |w|^{2^*(\sigma)}dx\biggl)^{2/{2^*(\sigma)}}}.
\end{equation}
Let $\eta \in \calC^\infty_c(F(B^N_{2r}))$ such that
$
0 \leq \eta \leq 1  \quad \textrm{ and } \quad \eta\equiv 1 \quad \textrm{in } \quad F(B^N_r).
$
Let
$$
w_n(p)=n^{\frac{N-2}{2}} w(n \rho(p))
$$
Then for $n  \in \N^*$ we define the test function defined in $M$ by
\begin{equation}
u_n(p)=\eta(p) w_n(p).
\end{equation}
We have by integration by parts
$$
\begin{array}{ll}
E(u_n):=&\displaystyle \int_\M |\n u_n|^2 dv_g-\l \int_\M u_n^2 dv_g\\\\
&\displaystyle= \int_{F(B^N_{2r})} \eta^2 |\n w_n|^2 dv_g-\int_{F(B^N_{2r})} (\eta \D \eta) w_n^2 dv_g-\l \int_{F(B^N_{2r})} u_n^2 dv_g \\\\
&\displaystyle \leq  \int_{F(B^N_{2r})}|\n w_n|^2 dv_g-\l \int_{F(B^N_{2r})} w_n^2 dv_g
+O \left( \int_{F(B^N_{2r}) \setminus F(B^N_r)} w_n^2 dv_g\right).
\end{array}
$$
Then by a change of variable formula and the fact that
$\rho(F(x))=|x|$
we obtain
\begin{equation}\label{Energy1}
E(u_n)\leq \int_{B^N_{2nr}} |\n w|^2 \sqrt{|g|}\left(\frac{x}{n}\right) dx -\frac{\l}{n^2} \int_{B^N_{2nr}} w^2 \sqrt{|g|}\left(\frac{x}{n}\right) dx+O \left(\frac{1}{n^2} \int_{B^N_{2nr}\setminus B^N_{nr}} w^2 dx\right).
\end{equation}
It's well known that, the components of the metric $g$ in the local chart of the exponential map are given by
$$
g_{ij}(x)=\delta_{ij}-\sum_{\a\b=1}^N \frac{R_{i\a j\b}(p_0)}{3}  x_\a x_\b+O(|x|^3) .
$$
Then Cartan expansion of the metric yields
\begin{equation}\label{DetMetric}
\sqrt{|g|}(x)=1-\frac{1}{6} \sum_{\a\b=1}^N R_{\a\b}(p_0) x_\a x_\b+O \left(|x|^3\right).
\end{equation}
Moreover for all $r_0$ a positive real and all $\a, \b=1,....,N$ we have
\begin{equation}\label{Mouhamed Paix et Salut sur Lui}
\int_{B^N_{r_0}} |\n w|^2 x_\a x_\b dx=\frac{\d_{\a\b}}{N} \int_{B^N_{r_0}} |x|^2 |\n w|^2 dx
\quad \textrm{and} \quad
\int_{B^N_{r_0}} w^2(x) x_\a x_\b dx=\frac{\d_{\a\b}}{N} \int_{B^N_{r_0}} |x|^2 w^2(x) dx.
\end{equation}
Then by \eqref{Energy1}, \eqref{DetMetric} and \eqref{Mouhamed Paix et Salut sur Lui}, we obtain
\begin{equation}\label{Mouhamed}
E(u_n)\leq \int_{\R^N} |\n w|^2 dx-\frac{S_g(p_0)}{6N n^2} \int_{B^N_{2nr}} |x|^2 |\n w|^2 dx -\frac{\l}{n^2} \int_{B^N_{2nr}} w^2 dx+O \left(\rho_1(n)\right)
\end{equation}
where
$$
\rho_1(n):=\frac{1}{n^3} \int_{B^N_{2nr}} |x|^3 |\n w|^2 dx+\frac{1}{n^2} \int_{B^N_{2nr}} w^2 dx+\frac{1}{n^4} \int_{B^N_{2nr}} |x|^2 w^2 dx
$$
Thanks to \eqref{GroundState}, it is easy follow that
\begin{equation}\label{Estimation Error}
\rho_1(n)=o\left(\frac{1}{n^2}\right) \quad \textrm{for} \quad N \geq 5 \qquad 
\textrm{ and } \qquad 
\rho_1(n)=O\left(\frac{1}{n^2}\right) \quad \textrm{for} \quad N =4.
\end{equation}
Morover
\begin{equation}\label{Error222}
\frac{1}{n^2} \int_{\R^N \setminus B^N_{2nr}} |x|^2 |\n w|^2 dx+\frac{1}{n^2} \int_{\R^N \setminus B^N_{2nr}} w^2 dx=o\left(\frac{1}{n^2}\right) \qquad \textrm{for}\quad N \geq 5.
\end{equation}
This implies that
\begin{equation}\label{EnergyEnergy}
\begin{cases}
\displaystyle E(u_n) \leq \int_{\R^N} |\n w|^2 dx-\frac{S_g(p_0)}{6N n^2} \int_{\R^N} |x|^2 |\n w|^2 dx -\frac{\l}{n^2} \int_{\R^N} w^2 dx+o\left(\frac{1}{n^2}\right)\qquad \qquad\textrm{for}\quad N \geq 5\\\\
\displaystyle E(u_n) \leq \int_{\R^N} |\n w|^2 dx-\frac{S_g(p_0)}{6N n^2} \int_{B^N_{2nr}} |x|^2 |\n w|^2 dx -\frac{\l}{n^2} \int_{B^N_{2nr}} w^2 dx+O \left(\frac{1}{n^2}\right) \qquad \textrm{ for}\quad N =4.
\end{cases}
\end{equation}
We have also that
$$
\int_\M \rho^{-\s} u_n^{2^*(\s)} dv_g=\int_{F(B^N_{nr})} \rho^{-\s} w_n^{2^*(\s)} dv_g+\int_{F(B^N_{2nr})\setminus F(B^N_{nr})} \rho^{-\s} u_n^{2^*(\s)} dv_g.
$$
By change of variable formula, we obtain
\begin{equation}\label{Denominateur}
\begin{array}{ll}
\displaystyle\int_\M \rho^{-\s} u_n^{2^*(\s)} dv_g&\displaystyle=\int_{B^N_{nr}} |x|^{-\s} w^{2^*(\s)} \sqrt{|g|}\left(\frac{x}{n}\right) dx+O \left(\int_{B^N_{2nr}\setminus B^N_{nr}} |x|^{-\s} w^{2^*(\s)} dx \right)\\\\
&\displaystyle =\int_{\R^N} |x|^{-\s} w^{2^*(\s)} dx-\frac{S_g(p_0)}{6N n^2} \int_{\R^N} |x|^{2-\s} w^{2^*(\s)} dx+O \left(\rho_2(n)\right)
\end{array}
\end{equation}
where 
$$
\rho_2(n)=\int_{B^N_{2nr}} |x|^{-\s} w^{2^*(\s)} dx+\frac{1}{n^2} \int_{\R^N\setminus B^N_{nr}} |x|^{2-\s} w^{2^*(\s)} dx+\frac{1}{n^3} \int_{B^N_{nr}} |x|^{3-\s} w^{2^*(\s)} dx.
$$
It's easy follows that for all $N \geq 4$,
$$
\rho_2(n)=o\left(\frac{1}{n^2}\right).
$$
Therefore by Taylor expansion we get that for $N \geq 4:$
\begin{equation}\label{DenomEstim}
\left(\int_\M \rho^{-\s} u_n^{2^*(\s)} dv_g \right)^{2/2^*(\s)}
=
\left(\int_{\R^N} |x|^{-\s} w^{2^*(\s)} dx\right)^{2/2^*(\s)}\left\lbrace 1-\frac{2}{2^*(\s)}\frac{S_g(p_0)}{6N n^2} \int_{\R^N} |x|^{2-\s} w^{2^*(\s)} dx+o\left(\frac{1}{n^2}\right)\right\rbrace.
\end{equation}
Hence from \eqref{Hardy-SobolevBestConstant}, \eqref{EnergyEnergy} and \eqref{DenomEstim}, we obtain
\begin{equation}
\mu_{\l,\s} \leq S_{N,\s}-\frac{S_g(p_0)}{6N n^2} \int_{\R^N} |x|^2 |\n w|^2 dx-\frac{\l}{n^2} \int_{\R^N} w^2 dx+\frac{2}{2^*(\s)} \frac{S_g(p_0)}{6N n^2} \int_{\R^N} |x|^{2-\s} w^{2^*(\s)} dx+o\left(\frac{1}{n^2}\right)
\end{equation}
Moreover the ground state solution $w$ of the Hardy-Sobolev best constant $S_{N,\s}$ satisfies
\begin{equation}\label{GSS}
-\D w=S_{N,\s} |x|^{-\s} w^{2^*(\s)-1} \qquad \quad \textrm{in } \R^N.
\end{equation}
Then we multiply the above equation by $|x|^2 w$ and we integrate by parts twice to get
$$
\int_{\R^N} |x|^2 |\n w|^2 dx=N \int_{\R^N} w^2 dx+S_{N,\s} \int_{\R^N} |x|^{2-\s} w^{2^*(\s)} dx.
$$
Using this with the fact that the parameter $\l$ is negative, we obtain
\begin{equation}\label{FinN5}
\mu_{\l,\s} \leq S_{N,\s}-\frac{S_g(p_0)+6\l}{6N n^2} \int_{\R^N} |x|^2 |\n w|^2 dx+o\left(\frac{1}{n^2}\right) \qquad \textrm{for } N \geq 5.
\end{equation}
For the case $N=4$, we let $\vp \in \calC^\infty_c(B_{3r}^N)$ such that 
$$
\vp \equiv 1 \qquad\textrm{in } B^N_{2r} \qquad\textrm{ and }\quad \D \vp, |\n \vp| \leq \textrm{Const}.
$$
Define $\vp_n(x)=\vp\left(\frac{x}{n}\right)$. We multiply \eqref{GSS} by
by $\vp_n |x|^2 w$ and we integrate by parts to get
$$
\begin{array}{ll}
\displaystyle\int_{B^N_{2nr}} |x|^2 |\n w|^2 dx-N\int_{B^N_{2nr}} w^2 dx=&\displaystyle S_{N,\s} \int_{B^N_{3nr}}\vp_n |x|^{2-\s} w^{2^*(\s)} dx\\\\
&\displaystyle+\frac{1}{2n^2} \int_{B^N_{3nr}\setminus B^N_{2nr}} w^2(|x|^2 \D \vp_n dx+ |x|(\n \vp_n \n |x|^2)) dx.
\end{array}
$$
By \eqref{GroundState}, we obtain the estimation
$$
\int_{B^N_{3nr}} |x|^{2-\s} w^{2^*(\s)} dx= \int_{\R^N}|x|^{2-\s} w^{2^*(\s)} dx-\int_{\R^N \setminus B^N_{3nr}}|x|^{2-\s} w^{2^*(\s)} dx=Const.+o(1).
$$
and
$$
\int_{B^N_{3nr}\setminus B^N_{nr}} w^2 dx=Const.+o(1).
$$
This yiels
$$
\int_{B^N_{2nr}} |x|^2 |\n w|^2 dx=N \int_{B^N_{2nr}} w^2 dx+o(1) \qquad \textrm{ for } N=4.
$$
Hence using this with \eqref{EnergyEnergy} and \eqref{DenomEstim}, we obtain
\begin{equation}\label{FinN4}
\mu_{\l,\s} \leq S_{N,\s}-\frac{S_g(p_0)+6\l}{6N n^2} \int_{B^N_{2nr}} |x|^2 |\n w|^2 dx+O\left(\frac{1}{n^2}\right) \qquad \textrm{for } N=4.
\end{equation}
Note that by \eqref{GroundState}, we have 
\begin{equation}\label{NdeyeSigaBeye}
\int_{\R^N} |x|^2 |\n w|^2 dx< \infty \quad\textrm{for } N \geq5 \qquad \textrm{and}\qquad \int_{B^N_{2nr}} |x|^2 |\n w|^2 dx=O\left(\textrm{log}(n)\right) \quad \textrm{ for } N=4.
\end{equation}
From \eqref{FinN5}, \eqref{FinN4} and \eqref{NdeyeSigaBeye} we obtain
\begin{equation}
\mu_{\l,\s}< S_{N,\s}
\quad \textrm{
provided that }\quad S_g(p_0) >-6\l.
\end{equation}
This ends the proof. \QED
\textbf{ACKNOWLEDGEMENTS:} I wish to thank my supervisor Mouhamed Moustapha Fall and professeur Frederic Robert for their help and usefull discussions. This work is supported by the German Academic Exchange Sercive(DAAD).\\\\

\end{document}